\documentclass[12pt]{article}
\usepackage{amssymb,latexsym, amsmath, amscd, array, graphicx}
\usepackage{amssymb}

\usepackage{setspace}
\input xy
\xyoption{all}
\makeatletter 

\usepackage{setspace}
\usepackage{enumitem}

\date{}

\newtheorem{theorem}{Theorem}[section]

\newtheorem{proposition}[theorem]{Proposition}

\newtheorem{problem}[theorem]{Problem}

\newcommand{\z}{{\Bbb Z}}

\newcommand{\re}{{\Bbb R}}

\newcommand{\invlim}{\raisebox{-1ex}{$\stackrel{\hbox{lim}}{\leftarrow}$}}

\newcommand{\lo}{\rightarrow}

\newcommand{\black}{{\blacksquare}}

\begin{document}

\title{\bf On Raymond -Williams' example}

\author{  Michael  Levin\footnote{This research was supported by 
THE ISRAEL SCIENCE FOUNDATION (grant No. 522/14) }}

\maketitle
\begin{abstract} Raymond and Wiliams  in  {\em  ``Examples of p-adic transformation groups''.} 
Ann. of Math. (2) 78 (1963) 92-106
constructed an action of  the $p$-adic integers $A_p$ on an $n$-dimensional compactum $X$, $n \geq 2$, with 
the orbit space $X/A_p$ of dimension $n+2$. 
We present a simpler construction of  such an example. 
\bigskip
\\
{\bf Math. Subj. Class.:}  55M10, 22C05 (54F45)
\end{abstract}
\begin{section}{Introduction} 
The Hilbert-Smith conjecture asserts that a compact group acting effectively on a manifold  must be a Lie group.
This conjecture is equivalent to the following  one: the group  $A_p$  of the $p$-adic integers cannot act effectively 
on a manifold. Yang \cite{yang} showed  that if $A_p$ acts effectively on a manifold $M$ then
$\dim M/A_p =\dim M +2$. In order to verify if the latter dimensional relation ever occurs in a more general setting,
Raymond and Williams \cite{raymond-williams} constructed an action of $A_p$  on an $n$-dimensional 
compactum (=compact metric space) $X, n \geq 2$, 
with $\dim X/A_p=n+2$. We present  a simpler approach for constructing such an example.  Let first note
that it is sufficient  to consider only the case $n=2$ because the induced action of $A_p$ on 
$X \times [0,1]^k$  defined by $g(x, t)=(gx, t), g \in A_p, x \in X$ and $t \in [0,1]^k$ 
provides the  corresponding  example with  $ \dim X \times [0,1]^k=k+2$ and 
$\dim (X \times [0,1]^k)/A_p =\dim(X/A_p)\times [0,1]^k=\dim X/A_p +k$. Thus without loss of generality 
we restrict ourselves to  the dimension  $n=2$ (although the construction below can be easily adjusted
to any finite dimension $n\geq 2$, it is slightly easier to visualize  it  for $n=2$).
Our example is constructed in the next section and a few related problems are posed in the last section.

\end{section}
\begin{section}{Example}

Consider  the unit $3$-sphere  $S^3$ in $\re^4$. Represent $\re^4$ as the product $\re^2 \times \re^2_\perp$
 of a coordinate plane $\re^2$ and its orthogonal complement $\re^2_\perp$ and let $S^1=S^3 \cap \re^2$ and
 $S^1_\perp =S^3 \cap \re^2_\perp$. Take the  closed $\epsilon$-neighborhood $F$ of $S^1$.
 Then $F$ and  $F_\perp=S^3 \setminus {\rm int} F$ can be represented as the products
 $F=S^1\times D$ of $S^1$ and 
  $F_\perp= D_\perp \times  S^1_\perp $
 with $D$ and $D_\perp$ being $2$-disks. Thus we have that 
 $\partial F = \partial F_\perp =S^1 \times \partial D = \partial D_\perp \times S^1_\perp$.
  Let $T=\re  /\z$ freely act on the circle $\partial D$ by rotations. Then this action induces  
the corresponding free action  on  $\partial F  = \partial F_\perp$ and the later action obviously extends  over 
$  F_\perp  $ as a free action and  it extends over $F$  by the  rotations of  $D$ induced by the rotations of $\partial D$.
This way we have defined an action of $T$ on $S^3$ whose  fixed point set  is $S^1$ and
$T$ acts freely on $S^3\setminus S^1$.

\begin{proposition}
\label{prop0}
Let $\Delta$ be a $4$-simplex. There is an unknotted circle $S^1$ in $S^3=\partial \Delta$ such that
for each $3$-simplex $\Delta'$ the intersection $S^1 \cap \Delta'$ is an interval whose end-points are  in 
$\partial \Delta'$ and whose  interior is in the interior of $\Delta'$.
\end{proposition}
{\bf Proof.} Fix a  $3$-simplex $\Delta'$ of $\Delta$.
First observe that $S^1$ can be embedded into $\partial \Delta'$ so that  $S^1$ does not
contain the vertices of $\Delta'$ and the intersection of $S^1$  with
every $2$-simplex $\Delta''$ of $\Delta'$ is an interval whose end-points are in
$\partial \Delta''$ and whose interior is in the interior of $\Delta''$. 
Take any two $3$-simplexes $\Delta'_-$ and $\Delta'_+$ different from $\Delta'$ and
take any   points  $x_- \in S^1 \cap \Delta'_- \cap \Delta'$ and  $x_+\in S^1 \cap \Delta'_+\cap \Delta'$  not lying
in $\partial (\Delta'_- \cap \Delta')$ and $\partial (\Delta'_+ \cap \Delta')$ respectively. Now for each $3$-simplex
$\Delta'_*$ of $\Delta$ different from $\Delta', \Delta'_-$ and $\Delta'_+$ push the interior of 
the interval $S^1 \cap \Delta'_*$ inside the interior $\Delta'_*$ not moving  the end-pints of the interval.
The interior of  the interval  of  $S^1$ connecting $x_-$ and $x_+$ inside $\Delta'_- \cup \Delta'_+$  push into
the interior of 
$\Delta'$ not moving the end-points $x_-$ and $x_+$. And finally push the intervals
$S^1\cap \Delta'_-$ and $S^1\cap \Delta'_+$ into the interior of 
$\Delta'_-$ and $\Delta'_+$ respectively not moving the end-points of the intervals.
Thus  $S^1 \subset \partial \Delta$  has the required properties. $\black$
\\
\\
From now we identify $S^3$ with the boundary of a $4$-smplex $\Delta$ and assume that  $S^1\subset S^3$
satisfy the conclusions of  Proposition \ref{prop0}.
\\
\\
 Recall that 
 $F$  is represented as  the product $F=S^1\times D$, take $5$ disjoint closed intervals 
 $I_1$, ..., $I_5$ in $S^1$  lying in the interior of different $3$-simplexes of $\Delta$
 and  denote  $B_i =I_i \times D'$ where $D' \subset D$ is a  concentric disk in $D$.
 Taking  $D'$ to be  small enough we may assume that
 $B_i$ lie in the interior of different $3$-simplexes of $\Delta$.
  Note that  each  $3$-ball 
 $B_i$ is invariant under the action of $T$
 and hence $M=S^3\setminus {\rm int}(B_1\cup\dots\cup B_5)$ is invariant as well. 

\begin{proposition}
\label{prop1}
${}$

(i) For every $g \in T$ the map $x \lo gx, x \in M,$  is ambiently isotopic
to a homeomorphism $\phi : M \lo M$ 
 by an isotopy of $M$  that does not move the points of $S^1\cup \partial  M$ 
and such that  $\phi$
restricted to the $2$-skeleton $\Delta^{(2)}$ coincides with
 the inclusion of $\Delta^{(2)}$ into $M$ and $\phi(\Delta' \cap M)=\Delta' \cap M$ for
 every $3$-simplex $\Delta'$  of $\Delta$.

(ii)  There is a retraction of $r: M\lo \Delta^{(2)}$ such that for every $3$-simplex $\Delta'$ of $\Delta$
we have  that $r(\Delta'  \cap M) \subset \partial \Delta'$.
\end{proposition}
{\bf Proof.} 

(i) Since $T$ is path-connected $g : M \lo M$ is isotopic to the identity map of $M$ by an isoptopy that does not
move the points of $M\cap S^1$ (the fixed point set of $T$). For each $\partial B_i \subset \partial M$ we can
change this isotopy on a neighborhood of $\partial B_i$ in $M$ to get in addition that the points of
$\partial B_i$ are not moved.

(ii) For  each $3$-simplex $\Delta'$ of $\Delta$  define $r$ on $\Delta' \cap M$
 as  the restriction of a radial projection to
$\partial \Delta'$ from a point 
of $\Delta'$ that
that does not belong to $M$.
$\black$
\\ \\
For a prime number $p$ 
consider the subgroup $\z_p=\z/p\z$  of $T=\re/\z$ and  let 
 $\Gamma_M =M/\z_p$.   Denote by $M_+$ the space which is the union of $p$ copies $M_a$, $a \in \z_p$, of $M$
 with $\partial M$ being identified in all the copies by the identity map. Thus
 all $M_a \subset M_+$ intersect each other at $\partial M$.  Define the action of $\z_p$ on $M_+$
 by sending  $x \in M_a$ to $gx \in M_{g+a}$.  Note that $\partial M \subset M_+$ is invariant under the action of
 $\z_p$ on $M_+$ and the natural  projection  $\alpha_+ : M_+ \lo M$ sending each $M_a$ to $M$ by the identity map
 is  equivariant.
 
 Let $\Gamma_+$ be the mapping cylinder of $\alpha_+$. The actions of $\z_p$ on $M_+$ and $M$
 induce the corresponding action of $\z_p$ on $\Gamma_+$. The spaces $M_+$ and $M$  can be considered
 as natural subsets of the mapping cylinder $\Gamma_+$,  we  will refer to $M_+$ and $M$ as the bottom 
 and the top of $\Gamma_+$  respectively. 
  
Let  $\Gamma=\Gamma_+ /\z_p$,  $\gamma_+ : \Gamma_+ \lo \Gamma$  the projection
  and 
   $\Gamma_{\partial M}$ and $\Gamma_M$  the images  in $\Gamma$
   of the bottom set $M_+$  and the top set $M$  of $\Gamma_+$  under  the map $\gamma_+$.
   Note that $\alpha_+: M_+ \lo M$ induces the corresponding map $\alpha : \Gamma_{\partial M} \lo \Gamma_M$
   for which  $\Gamma$ is the mapping cylinder of $\alpha$.
   Also note that $\Gamma_M = M/\z_p$ and  $\Gamma_{\partial M}$   is   the space obtained from $M$
   by collapsing to singletons  the orbits of the action of $\z_p$ on $M$  lying in $\partial M $.
   Thus   there is a natural projection $\mu : M \lo \Gamma_{\partial M}$ and  $ \Delta^{(2)} \subset M$
   can be considered as a subset of $\Gamma_{\partial M}$ as well (we can identify $M$ with one of the spaces 
   $M_a\subset M_+$ and  regard $\mu $ as the projection $M_+ \lo  \Gamma_{\partial M} =M_+/\z_p$
   restricted  to $M_a$).  
   
     Let us define  a map 
 $\delta: \Gamma \lo \Delta$ by sending  $\Gamma_{\partial M}$ to $\partial \Delta$  and
 $\Gamma_M$ to the barycenter of $\Delta$ such that $\delta$ is the identity map on $\Delta^{(2)}$,
$\delta(x) \in \Delta'$ for every $3$-simplex $\Delta'$ of $\Delta$ and $x \in \mu(M \cap \Delta')$
and $\delta$ is linearly extended along the intervals of the mapping cylinder $\Gamma$ (in our descriptions 
we abuse notations  regarding simplexes of $\Delta$ also as subsets  of other spaces involved). 
 Denote $\delta_+ =\delta \circ \gamma_+  : \Gamma_+ \lo \Delta$.

 In the notation below
   an element  $g \in \z_p$ that   appears
 in a composition with other maps is  regarded  as a homeomorphism of the space on which $\z_p$ acts.
 
\begin{proposition}
\label{prop2}
For every $3$-simplex $\Delta_i $ of $\Delta$  there is a map 
$r^i_+ :\delta_+^{-1}(\Delta_i) \lo \partial \Delta_i$ such that

(i)
 $r^i_+$ coincides with $\delta_+$  on  
 $ \delta_+^{-1}(\partial \Delta_i))$;

  (ii) for every triangulation of $\Gamma_+$ and  every  $g_1, \dots, g_5 \in \z_p$    there is a map
   $r_+  : \Gamma_+ \lo \partial \Delta$  such that  $r_+ (\Gamma_+^{(3)}) \subset \Delta^{(2)}$ and 
  $r_+$ restricted to $ \delta_+^{-1}(\Delta_i)$ coincides with $r^i_+ \circ g_i$
  for each $i$.

\end{proposition}
{\bf Proof.} 
Recall that $\Gamma_+$ is the mapping cylinder of $\alpha_+ : M_+ \lo M$ and let
$\beta_+ : \Gamma_+ \lo M$ be the projection  to the top set $M$ of $\Gamma_+$.
Also recall that $M_+$ is the union of $p$ copies $M_a$, $a \in \z_p$, of $M$  intersecting each other at $\partial M$.

Denote $\Delta^{(2)}_a =\delta_+^{-1}(\Delta^{(2)}) \cap M_a$,  
$(\Delta_i\cap M)_a =\delta_+^{-1}(\Delta_i)  \cap M_a$,
$(\partial \Delta_i)_a=\delta_+^{-1}(\partial \Delta_i)\cap M_a$
  and $(S^1\cap M)_a =\alpha_+^{-1}(S^1\cap M) \cap M_a$.
Note  that $(S^1\cap M)_a\cap (\Delta_i \cap M)_a$ splits into two intervals each of them connects in $M_a$ the sets
$\partial B_i$ and $(\partial \Delta_i)_a$. One of this intervals we will denote by $(S_i)_a$.

By   (i) of Proposition \ref{prop1} the map $\beta_+$ 
 can be isotoped  relative to $\beta_+^{-1}((S^1 \cap M)\cup \partial M)$ into a map 
 $\omega_+  : \Gamma_+ \lo M$  such that  $\omega_+$  and $\delta_+$ restricted
 to each $\Delta^{(2)}_a \subset M_+$ coincide.
 Consider the map $r$ from (ii) of Proposition \ref{prop1} and note that 
 $(r \circ \omega_+)(\Gamma_+) \subset \Delta^{(2)}$.
 Define the map $r_+^i$ as the map $r \circ \omega_+$ restricted   to $\delta_+^{-1}(\Delta_i)$.
\\
\\
(i) follows from (ii) of Proposition \ref{prop1}.
\\
\\
(ii)   
 For a subset $A \subset  M_+$ by the mapping cylinder 
of  $\alpha_+$ over $A$ we mean the mapping cylinder of $\alpha_+ : A \lo \alpha_+(A)$ which is
 a subset of $\Gamma_+$.
Let  $S \subset M_+$  be  the union of $(S_i)_a$ for all $i$ and $a$, and
  $B\subset M_+$ the union  of all the balls   $B_i \subset M_+$
Consider a CW-structure of $\Gamma_+$  for which the interiors of the $4$-cells are 
the interiors in $\Gamma_+$  of the mapping cylinders of
$\alpha_+$  over $(\Delta_i \cap M)_a\setminus (S_i)_a$ without the points belonging to $M_+$ and $M$. 
Thus the  $3$-skeleton of this CW-structure is the union of $M_+\cup M$ and  the mapping
cylinder of $\alpha_+$ over $\delta^{-1}(\Delta^{(2)})\cup S\cup B$.

Let be the map $\phi : M_+ \lo \partial \Delta$ be defined by the maps $r^i_+ \circ g_i$  on
each $\delta_+^{-1}(\Delta_i)$.
 Note that  the maps $\phi$   and $r \circ  \omega_+$  coincide on  
 $\delta^{-1}_+(\Delta^{(2)})\cup S$.
 Also note $\phi$ and  $r \circ \omega_+$ are homotopic on $B$
 by a homotopy relative to $B\cap S$. 
 Thus $r \circ \omega_+$ restricted to the union of $M \subset  \Gamma_+$ with 
 the mapping cylinder of $\alpha_+$ over
 $\delta_+^{-1}(\Delta^{(2)}) \cup B \cup S$
   can be homotoped
 to a map $\Phi$ such that $\phi$ and $\Phi$ restricted to  $\delta_+^{-1}(\Delta^{(2)}) \cup B \cup S$ coincide.
 Thus we can extend $\phi $ to a map $\phi_+$
  from  the $3$-skeleton of the CW-structure of $\Gamma_+$ to $\Delta^{(2)}$. 

  Then,
since $\Delta^{(2)}$   is contractible inside  $\partial \Delta$, we may extend $\phi_+$ to a map 
$r_+ : \Gamma_+ \lo \partial \Delta$. 
 For any triangulation  of $\Gamma_+$, the $3$-skeleton
of the  triangulation  can be pushed off the interiors of the $4$-cells of $\Gamma_+$ relative
to the $3$-skeleton of the CW-structure  of $\Gamma_+$
and  (ii) follows.
$\black$
\\
\\
Denote by $\Gamma_* $ the space obtained from $\Gamma_+$ by collapsing the fibers of $\gamma_+$
to singletons over the set $\delta^{-1}(\Delta^{(2)})$  and let the maps
$\gamma_* : \Gamma_* \lo \Gamma$, $\delta_* : \Gamma_* \lo \Delta$, 
$r^i_*: \delta^{-1}(\Delta_i)\lo \partial \Delta_i$ and 
 $r_* : \Gamma_* \lo \partial \Delta$ be induced by $\gamma_+, \delta_+, r^i_+$ and $r_+$ respectively and
 consider $\Gamma_*$ with the action of $\z_p$ induced by the action of $\z_p$  on
 $\Gamma_+$.  
 \begin{proposition}
 \label{prop*}
 ${}$

 (i)  The conclusions of 
 Proposition \ref{prop2} hold with  the subscript ``$+$'' 
   being replaced everywhere  by the subscript ``$*$''.
 
 (ii)  The fixed point set of the action of $\z_p$ on $\Gamma_*$ is $2$-dimensional and
  there are  triangulations of $\Gamma_*$ and $\Gamma$ and a subdivision 
 of $\Delta$ for which the action of $\z_p$ on $\Gamma_*$ is simplicial and the maps $\gamma_*$
  and $\delta$ are simplicial.
 \end{proposition}
 {\bf Proof.}
 (i) is obvious and (ii)
 can be derived from the  construction of the spaces and the maps involved. 
 $\black$
\\  \\
 Let $L$ be a finite  $4$-dimensional  simplicial complex and $\lambda  : L \lo \Delta$ a simplicial map
 such that $\lambda$ is $1$-to-$1$ on each simplex of  $L$. Denote by 
 $L'$ the pull-back space of the maps $\lambda$ and $\delta : \Gamma \lo \Delta$ and  by
 $\Omega: L' \lo L$ the pull-back of $\delta$.
 
 \begin{proposition}
 \label{prop3}
 The map $\Omega : L' \lo L$ induces an isomorphism  of $H_4(L'; \z_p)$ and
 $H_4(L; \z_p)$.
 \end{proposition}
 {\bf Proof.}  Recall that for every $3$-simplex $\Delta'$ of $\Delta$, $\delta^{-1}(\Delta')$
 is the mapping cylinder a  map of degree $p$ from $\partial \Delta'  $ to a $2$-sphere $S^2$ (the boundary of one of 
 the $3$-balls $B_i$).
 
 Also recall that  $\Gamma=\delta^{-1}(\Delta)$ is the mapping cylinder of the map $\alpha$ and from the definition of
 $\alpha$ one can also observe that 
 $H_4(\delta^{-1}(\Delta), \delta^{-1}(\partial \Delta); \z_p)=\z_p$
 and $\delta$ induces an isomorphism between 
 $H_4(\delta^{-1}(\Delta), \delta^{-1}(\partial \Delta); \z_p)$ and 
 $H_4(\Delta, \partial \Delta; \z_p)=\z_p$.
 
 And finally recall that $\delta$ is $1$-to-$1$ over the $2$-skeleton of $\Delta$. 
 Consider the long exact sequences of the pairs $(L', L'^{(3)})$  and $(L, L^{(3)})$ for the homology 
 with coefficients in $\z_p$. The facts above imply
 that $\Omega$ induces isomorphisms
 $H_3(L'^{(3)}; \z_p)\lo H_3(L^{(3)}; \z_p)$ and 
 $H_4(L', L'^{(3)}; \z_p)\lo H_4(L, L^{(3)}; \z_p)$. Then, by the $5$-lemma, $\Omega$
 induces an isomorphism $H_4(L'; \z_p)\lo H_4(L; \z_p)$ as well. $\black$
 \\
 
\begin{proposition}
\label{prop4}
 Let $G=\z_{p^k}$ act simplicially  on
a finite $4$-dimensional simplicial complex $K$ such that the action of $G$ is free on
$K \setminus K^{(2)}$.  Then there is a finite $4$-dimensional simplicial complex $K'$, a simplicial action
of $G'=\z_{p^{k+1}}$ on $K'$ and a map $\omega : K' \lo K$ such that the action of $G'$ is free on
$K'\setminus K'^{(2)}$ and

(i)  the actions of $G$ and $G'$  agree with $\omega$ and the natural  epimorphism  $h : G' \lo G$.
By this we mean that $\omega(g'x))=h(g')\omega(x)$ for for every $x \in K'$ and $g' \in G'$;

(ii) there is a map $\kappa  : K' \lo K^{(3)}$ such that
$\kappa(K'^{(3)}) \subset K^{(2)}$ and
 $\kappa(\omega^{-1}(\Delta_K))\subset \Delta_K$
 for every simplex $\Delta_K$ of $K$;
 
 (iii) the map $K'/G' \lo K/G$ determined by $\omega$ induces an isomorphism
 $H_4(K'/G'; \z_p) \lo H_4(K/G; \z_p)$.
 
 \end{proposition}
 {\bf Proof.} Replacing the triangulation of $K$ by a subdivision we may assume
 that $L=K/G$ is a simplicial complex, the projection  $\pi : K\lo L$ is a simplicial map and
 $L$ admits a simplicial map $\lambda_L : L \lo \Delta$ to a $4$-simplex $\Delta$ such that
 $\lambda_L$ is $1$-to-$1$ on each simplex of $L$.
 For every $4$-simplex  $\Delta_L$ of $L$  fix a $4$-simplex $\Delta_K$ of $K$ such that 
 $\pi(\Delta_K)=\Delta_L$ and denote by $K_-$ the union of the $3$-skeleton  $K^{(3)}$ of $K$
  with all the $4$-simplexes of $K$ that we fixed.
  Let $K'_-$ be  the pull-back space of the maps 
 $\lambda_L \circ \pi|K_- : K_- \lo \Delta$ and  $ \delta_*: \Gamma_* \lo \Delta$, 
 $\omega_- : K'_- \lo K$ the pull-back map of $\delta_*$ and
 $\lambda_*  : K'_- \lo \Gamma_*$ the pull back of $\lambda_L\circ \pi|K_-$.

Let  $g$ be a generator  of $G'$ and  $l : G'\lo \z_p$ an epimorphism.
 We will first  define the action of $G'$ on $\omega_-^{-1}(K^{(3)})$. 
 For each $3$-simplex $\Delta_L$ of $L $ define the action of $G'$ on $\omega_-^{-1}(\pi^{-1}(\Delta_L))$ as follows.
 Fix a   $3$-simplex $\Delta_K$ of $\pi^{-1}(\Delta_L)$ 
  and let  $x \in \omega_-^{-1}(\Delta_K)$. 
 Define $y=g^ t x$ for $ 1 \leq t \leq p^k-1$
 as the point $y \in K'$ such that $\omega_-(y) \in h(g^t)(\Delta_K)$ and 
 $\lambda_*(y)=\lambda_*(x)$, and for $t=p^k$ define $y=g^t x$ as
 the point $y \in \Delta_K$ such that $\lambda_*(y)=l(g)\lambda_*(x)$. 
 We do this  independently  for every $3$-simplex $\Delta_L$  of $L$ and this
 way define  the action of $g$ on $\omega_-^{-1}(K^{(3)})$. It is easy to see that  the action of $g$ is well-defined, 
 $g^t$ for $t={p^{k+1}}$  is the identity map of  $\omega_-^{-1}(K^{(3)})$
  and hence the action of $g$ defines the action  of $G'$ on  $\omega_-^{-1}(K^{(3)})$.
  Note that
 \\
  
(*)    for  $g^t \in G', t=p^k, $ and   $g_*=l(g)  \in \z_p$ 
we have  that 
 $\lambda_* \circ g^t$ and $g_* \circ \lambda_*$
 coincide on $\omega_-^{-1}(\partial \Delta_K)$ 
  for every $4$-simplex $\Delta_K$ of $K$.
 \\
 
  Now we  will enlarge $K'_-$ to  a space $K'$   and  extend the action of $G'$ over 
   $K'$.
  Let $\Delta_L$ be a $4$-simplex of $L$. Recall that we fixed
   a $4$-simplex $\Delta_K$ in  $\pi^{-1}(\Delta_L)$.  For every 
  $g'=g^t \in G', 1\leq p^k-1 $ attach  to $g'(\omega_-^{-1}(\partial \Delta_K))$  a copy of 
   the space $\omega_-^{-1}(\Delta_K)$
  (which is is in its turn a copy of   $\Gamma_*$) by identifying $g'(\omega_-^{-1}(\partial \Delta_K))$   with
  $\omega_-^{-1}(\partial \Delta_K)$  according to $g'$ and for $x \in \omega_-^{-1}(\Delta_K)$ define
  $g'x$ as the as the point corresponding to $x$ in the attached space. 
We will define the action of $g^t$, $t=p^k$, on $\omega_-^{-1}(\Delta_K)$ by 
$y=g^t x, x\in \omega_-^{-1}(\Delta_K)$ such that $y \in \omega_-^{-1}(\Delta_K)$  and $\lambda_* (y)=l(g)\lambda_*(x)$.
By (*) the action of $g$ on $\omega_-^{-1}(\Delta_K)$   agrees with  the action of $g$ on $\omega_-^{-1}(K^{(3)})$.
 We do the above procedure  independently for every  $4$-simplex $\Delta_L$ of $L$ and
this way we define the space $K'$ and the action of $G'$ on $K'$.  We extend $\omega_-$ and $\lambda_*$ to
the maps $\omega: K' \lo K$ and $\lambda'_K : K' \lo \Gamma$ by 
   $\omega(g^t x)=\omega_-(x)$  and $\lambda'_K (g^t x) =\gamma_*(\lambda_*(x))$  for $x$ in a fixed $4$-simplex 
   $\Delta_K$ and 
$1\leq t \leq p^k$.

 It is easy to verify that the action of $G'$ on $K'$ and the maps $\omega$ and $\lambda'_K$ are
well-defined and the conclusion (i) of the proposition holds. Moreover  $\lambda'_K\circ g'=\lambda'_K$ 
for every $g' \in G'$ and  hence $\lambda'_K$
 defines the corresponding map $\lambda'_L : L'=K'/G' \lo \Gamma$. Then $L'$ is 
the pull-back of the maps $\lambda_L : L \lo \Delta$ and $\delta  : \Gamma \lo \Delta$
with $\lambda'_L$ being the pull-back of $\lambda_L$ and the map $\Omega : L' \lo L$ induced by $\omega$ being
the pull-back of $\delta$.
Thus, by Proposition \ref{prop3}, the conclusion (iii) of the proposition  holds as well.

Consider any triangulation of $K'$  for which
  the preimages under $\omega$ of the simplexes of $K'$ are subcomplexes of $K'$. 
Then  the map 
$r_* : \Gamma_* \lo \Delta^{(3)}$ provided by  Propositions \ref{prop*}  and    \ref{prop2}
for $g_1=\dots=g_5=0 \in \z_p$ defines 
  the corresponding map 
  $\kappa_- : K'_- \lo K^{(3)}$ such that $\kappa(\omega_-^{-1}(K^{(3)})) \subset K^{(2)}$. 
  The construction above  and
    Propositions \ref{prop*}  and \ref{prop2} allow us to extend 
    the map $\kappa_-$ to a map $\kappa : K' \lo K^{(3)}$
   satisfying the conclusion (ii) of the proposition. Recall that  the  triangulation of $K$ 
    is a subdivision of 
the original triangulation of $K$.
    Replacing $\kappa$ by its composition 
  with the simplical approximation of the identity map of $K$ with respect to the new and original triangulations 
  of  $K$ we get that
   the conclusion  (ii) of the proposition holds.
  
The rest of the conclusions of the proposition follows from
(ii) of Proposition \ref{prop*}. $\black$
\\ \\
Now we are ready to construct our example. Set $K_0$  to be any finite $4$-dimensional simplicial complex with $H_4(K_0; \z_p)\neq 0$
and let the trivial group $\z_{p^0}=0$  trivially act on $K_0$.
Construct by induction on $i$ finite $4$-dimensional simplicial complexes $K_i$, an action of $\z_{p^i}$ on $K_i$ and maps 
$\omega_{i+1} : K_{i+1} \lo K_i$ so that $K_{i+1}$,
 the action of $\z_{p^{i+1}}$ on $K_{i+1}$ and $\omega_{i+1}$ satisfy the conclusions of Proposition \ref{prop4}
 with $K$, $K'$, $k$ and $\omega$ being replaced by $K_i$, $K_{i+1}$, $i$ and $\omega_{i+1}$ respectively.
 Consider $X=\invlim (K_i, \omega_i)$. Replacing   the triangulation of $K_i$
 by a sufficiently fine barycentric subdivision we can achieve, by (ii) of Proposition \ref{prop4},  that
 $\dim X \leq 2$.  The conclusions (i)  and (iii)  of Proposition \ref{prop4} imply that $A_p =\invlim Z_{p^i}$ acts on $X$ so that
 for $L_i=K_i/\z_{p^i}$ and the map $\Omega_{i+1} : L_{i+1} \lo L_i$ induced by $\omega_i$ we have that
 $Y=X/A_p=\invlim (L_i, \Omega_i)$  and $H^4(Y; \z_p)\neq 0$.  Thus $\dim Y =4$.  Note that from the construction
 of $K'$ in Proposition \ref{prop4} one can derive that $X$ contains the $2$-skeleton $K^{(2)}_i$ of each $K_i$ and
 hence $\dim X =2$.

\end{section}
\begin{section}{Problems}
It is very challenging to try to modify  the construction presented in this paper to approach  the following problems.
\begin{problem}
Does there exist an action of $A_p$ on a $1$-dimensional compactum with 
the $3$-dimensional orbit space?
\end{problem}

\begin{problem}
Does there exist an action of $A_p$ on an $n$-dimensional compactum $X$ with
an invariant subset $X' \subset X$ so that $\dim X' < n$,  $\dim X/A_p=n+2$ and
the action of $A_p$ is free on $X\setminus X'$?
\end{problem}

Note that the negative answer to  the last problem settles the Hilbert-Smith conjecture for actions with 
finite dimensional orbit spaces.

\end{section}

Michael Levin\\
Department of Mathematics\\
Ben Gurion University of the Negev\\
P.O.B. 653\\
Be'er Sheva 84105, ISRAEL  \\
 mlevine@math.bgu.ac.il\\\\
\end{document}